\documentclass[12pt]{amsart}
\begin{document}
\theoremstyle{plain}
\newtheorem{thm}{Theorem}[section]
\newtheorem*{thm*}{Theorem}
\newtheorem{prop}[thm]{Proposition}
\newtheorem*{prop*}{Proposition}
\newtheorem{lemma}[thm]{Lemma}
\newtheorem{cor}[thm]{Corollary}
\newtheorem*{conj*}{Conjecture}
\newtheorem*{cor*}{Corollary}
\newtheorem{defn}[thm]{Definition}
\theoremstyle{definition}
\newtheorem*{defn*}{Definition}
\newtheorem{rems}[thm]{Remarks}
\newtheorem*{rems*}{Remarks}
\newtheorem*{proof*}{Proof}
\newtheorem*{not*}{Notation}
\newcommand{\npartial}{\slash\!\!\!\partial}
\newcommand\Trig{\operatorname{Trig}}
\newcommand\maF{{\mathcal F}}
\newcommand\CAP{{\mathcal AP}}
\newcommand\ep{\epsilon}
\newcommand\cM{{\mathcal M}}
\newcommand\cnn{{\mathcal N}}
\newcommand{\IS}{\mathbf{S}}
\newcommand\cS{{\mathcal S}}
\newcommand\cA{{\mathcal A}}
\newcommand{\Heis}{\operatorname{Heis}}
\newcommand{\Solv}{\operatorname{Solv}}
\newcommand{\Spin}{\operatorname{Spin}}
\newcommand{\SO}{\operatorname{SO}}
\newcommand{\ind}{\operatorname{ind}}
\newcommand{\Index}{\operatorname{index}}
\newcommand{\ch}{\operatorname{ch}}
\newcommand{\rank}{\operatorname{rank}}
\newcommand\End{\operatorname{End}}
\newcommand{\abs}[1]{\lvert#1\rvert}
 \newcommand{\A}{{\mathcal A}}
        \newcommand{\D}{{\mathcal D}}
\newcommand{\HH}{{\mathcal H}}
\newcommand\wM{\widetilde{M}}
        \newcommand{\LL}{{\mathcal L}}
        \newcommand{\B}{{\mathcal B}}
        \newcommand{\K}{{\mathcal K}}
\newcommand{\oo}{{\mathcal O}}\newcommand\maH{{\mathcal H}}
         \newcommand{\PP}{{\mathcal P}}
        \newcommand{\s}{\sigma}
        \newcommand{\coker}{{\mbox coker}}
        \newcommand{\p}{\partial}
        \newcommand{\dd}{|\D|}
        \newcommand{\n}{\parallel}  
\newcommand{\bma}{\left(\begin{array}{cc}}
\newcommand{\ema}{\end{array}\right)}
\newcommand{\bca}{\left(\begin{array}{c}}
\newcommand{\eca}{\end{array}\right)}
\newcommand{\sr}{\stackrel}
\newcommand{\da}{\downarrow}
\newcommand{\tD}{\tilde{\D}}
        \newcommand{\R}{\mathbf R}
        \newcommand{\C}{\mathbf C}
        \newcommand{\h}{\mathbf H}
\newcommand{\Z}{\mathbf Z}
\newcommand{\N}{\mathbf N}
\newcommand{\tto}{\longrightarrow}
\newcommand{\ben}{\begin{displaymath}}
        \newcommand{\een}{\end{displaymath}}
\newcommand{\be}{\begin{equation}}
\newcommand{\ee}{\end{equation}}

        \newcommand{\bean}{\begin{eqnarray*}}
        \newcommand{\eean}{\end{eqnarray*}}
\newcommand{\nno}{\nonumber\\}
\newcommand{\bea}{\begin{eqnarray}}
        \newcommand{\eea}{\end{eqnarray}}
\newcommand{\supp}[1]{\operatorname{#1}}
\newcommand{\norm}[1]{\parallel\, #1\, \parallel}
\newcommand{\ip}[2]{\langle #1,#2\rangle}
\setlength{\parskip}{.3cm}
\newcommand{\nc}{\newcommand}
\nc{\nt}{\newtheorem}
\nc{\gf}[2]{\genfrac{}{}{0pt}{}{#1}{#2}}
\nc{\mb}[1]{{\mbox{$ #1 $}}}
\nc{\real}{{\mathbb R}}
\nc{\comp}{{\mathbb C}}
\nc{\ints}{{\mathbb Z}}
\nc{\Ltoo}{\mb{L^2({\mathbf H})}}
\nc{\rtoo}{\mb{{\mathbf R}^2}}
\nc{\slr}{{\mathbf {SL}}(2,\real)}
\nc{\slz}{{\mathbf {SL}}(2,\ints)}
\nc{\su}{{\mathbf {SU}}(1,1)}
\nc{\so}{{\mathbf {SO}}}
\nc{\hyp}{{\mathbb H}}
\nc{\disc}{{\mathbf D}}
\nc{\torus}{{\mathbb T}}
\newcommand{\tk}{\widetilde{K}}
\newcommand{\boe}{{\bf e}}\newcommand{\bt}{{\bf t}}
\newcommand{\vth}{\vartheta}
\newcommand{\CGh}{\widetilde{\CG}}
\newcommand{\db}{\overline{\partial}}
\newcommand{\tE}{\widetilde{E}}
\newcommand{\tr}{\mbox{tr}}
\newcommand{\ta}{\widetilde{\alpha}}
\newcommand{\tb}{\widetilde{\beta}}
\newcommand{\txi}{\widetilde{\xi}}
\newcommand{\hV}{\hat{V}}
\newcommand{\IC}{\mathbf{C}}
\newcommand{\IZ}{\mathbf{Z}}
\newcommand{\IP}{\mathbf{P}}
\newcommand{\IR}{\mathbf{R}}
\newcommand{\IH}{\mathbf{H}}
\newcommand{\IG}{\mathbf{G}}
\newcommand{\IN}{\mathbf{N}}
\newcommand{\CC}{{\mathcal C}}
\newcommand{\CD}{{\mathcal D}}
\newcommand{\CS}{{\mathcal S}}
\newcommand{\CG}{{\mathcal G}}
\newcommand{\CL}{{\mathcal L}}
\newcommand{\CO}{{\mathcal O}}
\nc{\ca}{{\mathcal A}}
\nc{\cag}{{{\mathcal A}^\Gamma}}
\nc{\cg}{{\mathcal G}}
\nc{\chh}{{\mathcal H}}
\nc{\ck}{{\mathcal B}}
\nc{\cl}{{\mathcal L}}
\nc{\cm}{{\mathcal M}}
\nc{\cn}{{\mathcal N}}
\nc{\NN}{{\mathcal N}}
\nc{\cs}{{\mathcal S}}
\nc{\cz}{{\mathcal Z}}
\nc{\sind}{\sigma{\rm -ind}}
\newcommand{\la}{\langle}
\newcommand{\ra}{\rangle}
\newcommand{\lp}{L_p({\mathcal N},\tau )}
\newcommand{\vlp}{\Vert _{_{L_p({\mathcal N},\tau )}}}
\newcommand{\al}{\alpha}
\def\Aut{\operatorname{Aut}}
\def\T{\mathbb T}

\def\cross#1{\rlap{\hskip#1pt\hbox{$-$}}}
        \def\intcross{\cross{0.3}\int}
        \def\bigintcross{\cross{2.3}\int}

 \title{Semifinite Spectral Triples Associated with Graph $C^*$-Algebras} 

 \vspace{.5 in}

\author{}
\begin{abstract}
We review the recent construction of semifinite spectral triples for
graph $C^*$-algebras. These examples have inspired many other
developments and we review some of these such as the relation between
the semifinite index and the Kasparov product, examples of
noncommutative manifolds, and an index theorem in twisted cyclic theory
using a KMS state.
\end{abstract}
\maketitle
\begin{center}
{\bf Alan L. Carey}\\Mathematical Sciences Institute\\
Australian National University\\
Canberra, ACT. 0200, AUSTRALIA\\
e-mail: acarey@maths.anu.edu.au

{\bf John Phillips}\\
Department of Mathematics and Statistics,\\
University of Victoria\\Victoria, B.C. V8W 3P4, CANADA

{\bf Adam Rennie}\\
Institute for Mathematical Sciences\\University of Copenhagen\\
Universitetsparken 5, DK-2100, Copenhagen, DENMARK

All authors were supported by the Australian Research Council (Australia)
while the second   author was supported by NSERC (Canada), and the
third author by SNF (Denmark).
\end{center}     
\allowdisplaybreaks

\tableofcontents

\section{Introduction}

The extension of the Connes-Moscovici local index theorem in noncommutative
geometry to the case
of semifinite spectral triples (in which the bounded operators
on a Hilbert space are replaced by a general semifinite von Neumann 
algebra) in \cite{CPRS2} and \cite{CPRS3} leads naturally to the question of
finding interesting novel examples. 
In this article we will 
direct our attention to  some recent investigations of new noncommutative
semifinite spectral triples arising from graph $C^*$-algebras. Not all graph $C^*$-algebras
lead to spectral triples and we describe in some detail
by the example of the Cuntz algebra how badly things fail. The results we 
describe for the Cuntz algebra are by way of an announcement of
a more extensive program in which we study KMS-Dixmier functionals
and twisted cyclic theory to shed light on situations where there is no
faithful trace on the algebra. 

We attempt to give enough detail on notations and definitions to make the 
account reasonably self contained. We do not try to give proofs referring
to the original articles or preprints instead, \cite{CNNR,CPR1,CPR2,CRT}. Thus the 
second half of the
Introduction will be devoted to the definition of a semifinite spectral triple and some of the analytic subtleties of which one needs to take
account in this more general situation. We then describe the 
semifinite version of the Connes-Moscovici local index theorem in Section 2.
A brief review of graph $C^*$-algebras 
occupies Section 3 and their Kasparov modules
then follow in Section 4. 
The construction of semifinite spectral triples 
for graph algebras and the index pairing via spectral flow is
described in Section 5 (these Sections follow \cite{pr}).

We then describe some results which are the subject of several papers 
in preparation. The first of these follows \cite{KNR}
where we explain the relationship of semifinite spectral triples to
KK theory (Section 6).
The explanation of the verification of the axioms 
of \cite{C1}
for a noncommutative spin manifold for 
the case of certain graph $C^*$-algebras
along the lines described in \cite{PRen} is next in Section 7.
For algebras that do not posses a faithful trace a completely different
picture emerges. We choose to illustrate this using the Cuntz algebra
in Section 8. The 
first new ingredient is that, imitating the procedure
in \cite{pr}, 
using the Kasparov module leads to a non-summable spectral triple.
By modifying the trace to make it a weight we restore summability
but at the expense of having to study twisted cyclic cocycles.
Moreover, we find in Section 9 that, there is a natural replacement for
$K_1$ of the Cuntz algebra which we term `modular
$K_1$'. Modular $K_1$ pairs, via the 
semifinite local index theorem, with these twisted cocycles.

{\bf Acknowledgements} It is a pleasure to acknowledge the advice, work and
scholarship of our colleagues Jens Kaad, Nigel Higson, Ryszard Nest, David Pask,
Aidan Sims, Fyodor Sukochev.

\subsection{Semifinite spectral triples}

We begin with some semifinite versions of standard definitions and
results. Let $\tau$ be a fixed faithful, normal, semifinite trace
on a von Neumann algebra ${\mathcal N}$. Let ${\mathcal
K}_{\mathcal N }$ be the $\tau$-compact operators in ${\mathcal
N}$ (that is the norm closed ideal generated by the projections
$E\in\mathcal N$ with $\tau(E)<\infty$).

\begin{defn} A semifinite
spectral triple $(\A,\HH,\D)$ is given by a Hilbert space $\HH$, a
$*$-algebra $\A\subset \cn$ where $\cn$ is a semifinite von
Neumann algebra acting on $\HH$, and a densely defined unbounded
self-adjoint operator $\D$ affiliated to $\cn$ such that

1) $[\D,a]$ is densely defined and extends to a bounded operator
in $\cn$ for all $a\in\A$

2) $(\lambda-\D)^{-1}\in\K_\cn$ for all $\lambda\not\in{\R}$

3) The triple is said to be even if there is $\Gamma\in\cn$ such
that $\Gamma^*=\Gamma$, $\Gamma^2=1$,  $a\Gamma=\Gamma a$ for all
$a\in\A$ and $\D\Gamma+\Gamma\D=0$. Otherwise it is odd.
\end{defn}

\begin{defn}\label{qck} A semifinite spectral triple
$(\A,\HH,\D)$ is $QC^k$ for $k\geq 1$
($Q$ for quantum) if for all $a\in\A$ the operators $a$ and
$[\D,a]$ are in the domain of $\delta^k$, where
$\delta(T)=[\dd,T]$ is the partial derivation on $\cn$ defined by
$\dd$. We say that $(\A,\HH,\D)$ is $QC^\infty$ if it is $QC^k$
for all $k\geq 1$.
\end{defn}

\noindent{\bf Remarks concerning derivations and commutators}.  By
partial derivation we mean that $\delta$ is defined on some
subalgebra of $\cn$ which need not be (weakly) dense in $\cn$.
More precisely, $\mbox{dom}\ \delta=\{T\in\cn:\delta(T)\mbox{ is
bounded}\}$. 
One of the difficulties
we face in the semifinite case is that zero may be in the 
continuous spectrum of $\D$ and for that reason
we have to work with $(1+\D^2)^{1/2}$
instead. Then if $T\in{\mathcal N}$, one can show
that $[\dd,T]$ is bounded if and only if $[(1+\D^2)^{1/2},T]$ is
bounded, by using the functional calculus to show that
$\dd-(1+\D^2)^{1/2}$ extends to a bounded operator in $\cn$. In
fact, writing $\dd_1=(1+\D^2)^{1/2}$ and $\delta_1(T)=[\dd_1,T]$
we have \ben \mbox{dom}\ \delta^n=\mbox{dom}\ \delta_1^n\ \ \ \
\mbox{for all}\ n.\een We also observe that if $T\in\cn$ and
$[\D,T]$ is bounded, then $[\D,T]\in\cn$. Similar comments apply
to $[\dd,T]$, $[(1+\D^2)^{1/2},T]$. The proofs of these statements
can be found in \cite{CPRS2}.
\begin{defn} A $*$-algebra $\A$ is smooth if it is Fr\'{e}chet
and $*$-isomorphic to a proper dense subalgebra $i(\A)$ of a
$C^*$-algebra $A$ which is stable under the holomorphic functional
calculus.\end{defn}

Thus saying that $\A$ is \emph{smooth} means that $\A$ is
Fr\'{e}chet and a pre-$C^*$-algebra. Asking for $i(\A)$ to be a
{\it proper} dense subalgebra of $A$ immediately implies that the
Fr\'{e}chet topology of $\A$ is finer than the $C^*$-topology of
$A$ (since Fr\'{e}chet means locally convex, metrizable and
complete.) We will sometimes speak of $\overline{\A}=A$,
particularly when $\A$ is represented on Hilbert space and the
norm closure $\overline{\A}$ is unambiguous. At other times we
regard $i:\A\hookrightarrow A$ as an embedding of $\A$ in a
$C^*$-algebra. We will use both points of view.

It has been shown that if $\A$ is smooth in $A$ then $M_n(\A)$ is
smooth in $M_n(A)$, \cite{GVF,LBS}. This ensures that the
$K$-theories of the two algebras are isomorphic, the isomorphism
being induced by the inclusion map $i$. 
Moreover a smooth algebra has a sensible spectral theory which agrees with that
defined using the $C^*$-closure, and the group of invertibles is
open.

\begin{lemma}[\cite{R1}]\label{smo} If $(\A,\HH,\D)$ is a $QC^\infty$ spectral triple, then
$(\A_\delta,\HH,\D)$ is also a $QC^\infty$ spectral triple, where
$\A_\delta$ is the completion of $\A$ in the locally convex
topology determined by the seminorms \ben
q_{n,i}(a)=\n\delta^nd^i(a)\n,\ \ n\geq 0,\ i=0,1,\een where
$d(a)=[\D,a]$. Moreover, $\A_\delta$ is a smooth algebra.
\end{lemma}

We call the topology on $\A$ determined by the seminorms $q_{ni}$
of Lemma \ref{smo} the $\delta$-topology. Thus whenever we have a
$QC^\infty$ spectral triple (semifinite or not) it extends to a smooth
algebra. This is a necessary ingredient in order to define topological
cyclic homology.


\section{Summability and the Semifinite Local Index Theorem}
In the following, let $\mathcal N$ be a semifinite von Neumann
algebra with faithful normal trace $\tau$. Recall from \cite{FK}
that if $S\in\mathcal N$, the \emph{ t-th generalized singular
value} of S for each real $t>0$ is given by
$$\mu_t(S)=\inf\{||SE||\ : \ E \mbox{ is a projection in }
{\mathcal N} \mbox { with } \tau(1-E)\leq t\}.$$

The ideal $\LL^1({\mathcal N},\tau)$ consists of those operators $T\in
{\mathcal N}$ such that $\n T\n_1:=\tau( |T|)<\infty$ where
$|T|=\sqrt{T^*T}$. In the Type I setting this is the usual trace
class ideal. We will denote the norm on $\LL^1(\cn,\tau)$ by
$\n\cdot\n_1$. An alternative definition in terms of singular
values is that $T\in\LL^1(\cn,\tau)$ if $\|T\|_1:=\int_0^\infty \mu_t(T) dt
<\infty.$

When ${\mathcal N}\neq{\mathcal
B}({\mathcal H})$, $\LL^1(\cn,\tau)$ need not be complete in this norm but it is
complete in the norm $||.||_1 + ||.||_\infty$. (where
$||.||_\infty$ is the uniform norm). Another important ideal for
us is the domain of Dixmier traces:
$${\mathcal L}^{(1,\infty)}({\mathcal N},\tau)=
\left\{T\in{\mathcal N}\ : \Vert T\Vert_{_{{\mathcal
L}^{(1,\infty)}}} :=   \sup_{t> 0}
\frac{1}{\log(1+t)}\int_0^t\mu_s(T)ds<\infty\right\}.$$

 The reader
should note that ${\mathcal L}^{(1,\infty)}(\cn,\tau)$ is often taken to
mean an ideal in the algebra $\widetilde{\mathcal N}$ of
$\tau$-measurable operators affiliated to ${\mathcal N}$. Our
notation is however consistent with that of \cite{C} in the
special case ${\mathcal N}={\mathcal B}({\mathcal H})$. With this
convention the ideal of $\tau$-compact operators, ${\mathcal
  K}({\mathcal N})$,
consists of those $T\in{\mathcal N}$ (as opposed to
$\widetilde{\mathcal N}$) such that \ben \mu_\infty(T):=\lim
_{t\to \infty}\mu_t(T)  = 0.\een

\begin{defn}\label{summable} A semifinite  spectral triple
$(\A,\HH,\D)$ relative to $(\cn,\tau)$
with $\A$ unital is
$(1,\infty)$-summable if \ben
(\D-\lambda)^{-1}\in\LL^{(1,\infty)}(\cn,\tau)\ \ \ \mbox{for all}\ \ \
\lambda\in\C\setminus\R.\een
\end{defn}

We need to briefly discuss the notion of Dixmier trace. For
more information on semifinite Dixmier traces, see \cite{CPS2}.
For $T\in\LL^{(1,\infty)}(\cn,\tau)$, $T\geq 0$, the function \ben
F_T:t\to\frac{1}{\log(1+t)}\int_0^t\mu_s(T)ds \een is bounded. For
certain $\omega\in L^\infty(\R_*^+)^*$, \cite{CPS2,C}, we
obtain a trace on $\LL^{(1,\infty)}(\cn,\tau)$ by setting
$$ \tau_\omega(T)=\omega(F_T),\ \  T\geq 0$$ 
and extending to $\LL^{(1,\infty)}(\cn,\tau)$
by linearity.
We will not go into the properties that $\omega$ must satisfy
to give a trace leaving that to the well established literature see 
\cite{CPS2}.
However,
 given such an $\omega$ associated to the semifinite normal trace $\tau$,
we denote by $\tau_\omega$ the corresponding Dixmier trace.

The Dixmier trace $\tau_\omega$ vanishes on the ideal of trace class
operators. Moreover whenever the function $F_T$ has a limit at infinity,
all Dixmier traces return the limit as their value. 
This leads to the notion of {\it measurable operator} (see \cite{C}
for the type I case). We
refer to \cite{LSS} for the definitive discussion of
the notion of measurable operator 
in the semifinite setting.
All Dixmier traces take the same value on measurable operators 
and we use the notation 
$\bigintcross$ in this case. 
That is, if $T\in\LL^{(1,\infty)}(\cn,\tau)$ is measurable, for
any allowed functional $\omega\in L^\infty(\R_*^+)^*$ we have
$$\tau_\omega(T)=\omega(F_T)=\bigintcross T.$$

We now introduce (a special case of) 
the analytic spectral flow formula of 
\cite{CP1,CP2}. 
This formula starts with a semifinite spectral triple $(\A,\HH,\D)$
and computes the spectral flow from $\D$ to $u\D u^*$, 
where $u\in\A$ is unitary with $[\D,u]$ bounded, in the case
where $(\A,\HH,\D)$ is of dimension $p\geq 1$.
Thus for any $n>p$ we have by Theorem 9.3 of \cite{CP2}:
\be sf(\D,u\D 
u^*)=\frac{1}{C_{n/2}}\int_0^1\tau(u[\D,u^*](1+(\D+tu[\D,u^*])^2)^{-n/2})dt,
\label{basicformula}\ee
with $C_{n/2}=\int_{-\infty}^\infty(1+x^2)^{-n/2}dx$. 
This real number $sf(\D,u\D u^*)$
recovers the pairing of the $K$-homology class
$[\D]$ of $\mathcal A$ with the $K_1({\mathcal A})$ class $[u]$
\cite{CPRS2}. There is a geometric way to view this formula. It is shown in 
\cite {CP2} that
the functional $X\mapsto \tau(X(1+(\D+X)^2)^{-n/2})$ on ${\mathcal N}_{sa}$
determines an exact one-form
on an affine space modelled on  ${\mathcal N}_{sa}$. Thus (\ref{basicformula})
represents the integral of this one-form along the path 
$\{\D_t=(1-t)\D+ tu\D u^*\}$
provided one appreciates that $\dot\D_t=u[\D,u^*]$ is a tangent vector
to this path.

Next we remind the reader of the local index theorem in noncommutative geometry
or at least the special case needed for this paper.
The original type I version of this result is due to Connes-Moscovici
\cite{CM}. There are two new proofs, one due to Higson \cite{Hig} and one
in \cite{CPRS2}. The latter argument handles the case of semifinite
spectral triples. 
 In the
simplest terms, the local index theorem provides a formula for the
pairing of a finitely summable spectral triple $(\A,\HH,\D)$ with
the $K$-theory of $\overline{\A}$. The special case that we
require in this paper is as follows. 
\begin{thm}[\cite{CPRS2}] Let $(\A,\HH,\D)$ be an odd $QC^\infty$
$(1,\infty)$-summable semifinite spectral triple, relative
to $(\cn,\tau)$. Then for $u\in\A$ unitary the pairing of $[u]\in
K_1(\overline{\A})$ with $(\A,\HH,\D)$ is given by
$$ \la
[u],(\A,\HH,\D)\ra=\lim_{r\to 0^+}
r\tau(u[\D,u^*](1+\D^2)^{-1/2-r}).$$ In particular, the
limit on the right exists.
\end{thm}
For more information on this result  see \cite{CPS2,CPRS2,CPRS3,CM}.

\section{The Gauge Spectral Triple of a Graph $C^*$-Algebra}

 For a more detailed introduction to graph
$C^*$-algebras we refer the reader to \cite{BPRS,kpr,Rae} and the
references therein. A directed graph $E=(E^0,E^1,r,s)$ consists of
countable sets $E^0$ of vertices and $E^1$ of edges, and maps
$r,s:E^1\to E^0$ identifying the range and source of each edge.
{\bf We will always assume that the graph is} {\bf row-finite}
which means that each vertex emits at most finitely many edges.
Later we will also assume that the graph is \emph{locally finite}
which means  it is row-finite and each vertex receives at most
finitely many edges. We write $E^n$ for the set of paths
$\mu=\mu_1\mu_2\cdots\mu_n$ of length $|\mu|:=n$; that is,
sequences of edges $\mu_i$ such that $r(\mu_i)=s(\mu_{i+1})$ for
$1\leq i<n$.  The maps $r,s$ extend to $E^*:=\bigcup_{n\ge 0} E^n$
in an obvious way. A \emph{loop} in $E$ is a path $L \in E^*$ with
$s ( L ) = r ( L )$, we say that a loop $L$ has an exit if there
is $v = s ( L_i )$ for some $i$ which emits more than one edge. If
$V \subseteq E^0$ then we write $V \ge w$ if there is a path $\mu
\in E^*$ with $s ( \mu ) \in V$ and $r ( \mu ) = w$
(we also sometimes say
 that $w$ is downstream from $V$). A \emph{sink}
is a vertex $v \in E^0$ with $s^{-1} (v) = \emptyset$, a
\emph{source} is a vertex $w \in E^0$ with $r^{-1} (w) =
\emptyset$.

A \emph{Cuntz-Krieger $E$-family} in a $C^*$-algebra $B$ consists
of mutually orthogonal projections $\{p_v:v\in E^0\}$ and partial
isometries $\{S_e:e\in E^1\}$ satisfying the \emph{Cuntz-Krieger
relations}
\begin{equation*}
S_e^* S_e=p_{r(e)} \mbox{ for $e\in E^1$} \ \mbox{ and }\
p_v=\sum_{\{ e : s(e)=v\}} S_e S_e^*   \mbox{ whenever $v$ is not
a sink.}
\end{equation*}

It is proved in \cite[Theorem 1.2]{kpr} that there is a universal
$C^*$-algebra $C^*(E)$ generated by a non-zero Cuntz-Krieger
$E$-family $\{S_e,p_v\}$.  A product
$S_\mu:=S_{\mu_1}S_{\mu_2}\dots S_{\mu_n}$ is non-zero precisely
when $\mu=\mu_1\mu_2\cdots\mu_n$ is a path in $E^n$. Since the
Cuntz-Krieger relations imply that the projections $S_eS_e^*$ are
also mutually orthogonal, we have $S_e^*S_f=0$ unless $e=f$, and
words in $\{S_e,S_f^*\}$ collapse to products of the form $S_\mu
S_\nu^*$ for $\mu,\nu\in E^*$ satisfying $r(\mu)=r(\nu)$ (cf.\
\cite[Lemma
  1.1]{kpr}).
Indeed, because the family $\{S_\mu S_\nu^*\}$ is closed under
multiplication and involution, we have
\begin{equation}
C^*(E)=\overline{\mbox{span}}\{S_\mu S_\nu^*:\mu,\nu\in E^*\mbox{ and
}r(\mu)=r(\nu)\}.\label{spanningset}
\end{equation}
The algebraic relations and the density of $\mbox{span}\{S_\mu
S_\nu^*\}$ in $C^*(E)$ play a critical role throughout the paper.
We adopt the conventions that vertices are paths of length 0, that
$S_v:=p_v$ for $v\in E^0$, and that all paths $\mu,\nu$ appearing
in (\ref{spanningset}) are non-empty; we recover $S_\mu$, for
example, by taking $v=r(\mu)$, so that $S_\mu S_v^*=S_\mu
p_{r(\mu)}=S_\mu$.

If $z\in S^1$, then the family $\{zS_e,p_v\}$ is another
Cuntz-Krieger $E$-family which generates $C^*(E)$, and the
universal property gives a homomorphism $\gamma_z:C^*(E)\to
C^*(E)$ such that $\gamma_z(S_e)=zS_e$ and $\gamma_z(p_v)=p_v$.
The homomorphism $\gamma_{\overline z}$ is an inverse for
$\gamma_z$, so $\gamma_z\in\Aut C^*(E)$, and a routine
$\epsilon/3$ argument using (\ref{spanningset}) shows that
$\gamma$ is a strongly continuous action of $S^1$ on $C^*(E)$. It
is called the \emph{gauge action}. Because $S^1$ is compact,
averaging over $\gamma$ with respect to normalised Haar measure
gives an expectation $\Phi$ of $C^*(E)$ onto the fixed-point
algebra $C^*(E)^\gamma$:
\[
\Phi(a):=\frac{1}{2\pi}\int_{S^1} \gamma_z(a)\,d\theta\ \mbox{ for
}\ a\in C^*(E),\ \ z=e^{i\theta}.
\]
The map $\Phi$ is positive, has norm $1$, and is faithful in the
sense that $\Phi(a^*a)=0$ implies $a=0$.

From Equation (\ref{spanningset}), it is easy to see that a graph
$C^*$-algebra is unital if and only if the underlying graph is
finite. When we consider infinite graphs, formulas which involve
sums of projections may contain infinite sums. To interpret these,
we use strict convergence in the multiplier algebra of $C^*(E)$:

\begin{lemma}[\cite{Rae}]\label{strict}
Let $E$ be a row-finite graph, let $A$ be a $C^*$-algebra
generated by a Cuntz-Krieger $E$-family $\{T_e,q_v\}$, and let
$\{p_n\}$ be a sequence of projections  in $A$. If $p_nT_\mu
T_\nu^*$ converges for every $\mu,\nu\in E^*$, then $\{p_n\}$
converges strictly to a projection $p\in M(A)$.
\end{lemma}

\section{Kasparov modules for graph algebras}

The Kasparov modules considered here are for $C^*$-algebras with trivial grading. 
\begin{defn} An odd Kasparov $A$-$B$-module consists of a
countably generated ungraded right $B$-$C^*$-module, with $\phi:A\to
End_B(E)$ a $*$-homomorphism, together with $P\in End_B(E)$  such that
$a(P-P^*),\ a(P^2-P),\ [P,a]$ are all compact
endomorphisms. Alternatively we may take $V=2P-1$ in favour of $P$
such that $a(V-V^*),\ a(V^2-1),\ [V,a]$ are all compact endomorphisms for
all $a\in A$.
\end{defn}

By
\cite[Lemma 2, Section 7]{K}, the pair $(\phi,P)$ determines a
$KK^1(A,B)$ class, and every class has such a representative. The
equivalence relations on pairs $(\phi,P)$ that give $KK^1$ classes
are unitary equivalence $(\phi,P)\sim (U\phi U^*,UPU^*)$ and
homology: $(\phi_1,P_1)\sim (\phi_2,P_2)$ if $P_1\phi_1(a)-P_2\phi_2(a)$ is a
compact endomorphism for all $a\in A$. The latter may be recast in
terms of operator homotopies as well, see \cite[Section 7]{K}.

We recall the construction of an odd Kasparov module for certain graph
$C^*$-algebras from \cite{pr}.
For $E$ a row finite directed graph, we set $A=C^*(E)$,
$F=C^*(E)^\gamma$, the fixed point algebra for the $U(1)$ gauge
action. The algebras $A_c,F_c$ are defined as the finite linear
span of the generators.

Right multiplication makes $A$ into a right $F$-module, and
similarly $A_c$ is a right module over $F_c$. We define an
$F$-valued inner product $(\cdot|\cdot)_R$ on both these modules
by
$$ (a|b)_R:=\Phi(a^*b),$$
where $\Phi$ is the canonical expectation $A\to F$.
\begin{defn}\label{mod} Let $X$ be the right $F$ $C^*$-module obtained by
completing $A$ (or $A_c$) in the norm
$$ \Vert x\Vert^2_X:=\Vert (x|x)_R\Vert_F=\Vert
\Phi(x^*x)\Vert_F.$$
\end{defn}

The algebra $A$ acting by multiplication on the left of $X$
provides a representation of $A$ as adjointable operators on $X$.
We let $X_c$ be the copy of $A_c\subset X$. The $\T^1$ action on $X_c$ is
unitary and extends to a strongly continuous unitary action on $X$.

For each $k\in\Z$, the spectral projection onto the $k$-th spectral subspace
of the $\T^1$ action is the operator $\Phi_k$ on $X$  by
$$\Phi_k(x)=\frac{1}{2\pi}\int_{\T^1}z^{-k}\gamma_z(x)d\theta,\ \ z=e^{i\theta},\ \ x\in X.$$
Observe that on generators we have \begin{equation}\Phi_k(S_\al
S_\beta^*)=\left\{\begin{array}{lr}S_\al S_\beta^* & \ \
|\al|-|\beta|=k\\0 & \ \ |\al|-|\beta|\neq
k\end{array}\right..\label{kthproj}\end{equation}

We quote the following result from \cite{pr}.
\begin{lemma}\label{phiendo} The operators $\Phi_k$ are adjointable endomorphisms
of the $F$-module $X$ such that $\Phi_k^*=\Phi_k=\Phi_k^2$ and
$\Phi_k\Phi_l=\delta_{k,l}\Phi_k$. If $K\subset\Z$ then the sum
$\sum_{k\in K}\Phi_k$ converges strictly to a projection in the
endomorphism algebra. The sum $\sum_{k\in\Z}\Phi_k$ converges to
the identity operator on $X$.
\end{lemma}

The unbounded operator of the next proposition is the generator of the
$\T^1$ action on $X$.
We refer to Lance's book, \cite[Chapters 9,10]{L}, for information
on unbounded operators on $C^*$-modules.

\begin{prop}\label{dee}[cf \cite{pr}] Let $X$ be the right $C^*$-$F$-module of
Definition \ref{mod}. Define $X_\D\subset X$ to be the linear space
$$ X_\D=
\{x=\sum_{k\in\Z}x_k\in X:\Vert \sum_{k\in\Z}k^2(x_k|x_k)_R\Vert<\infty\}.$$
For $x\in X_\D$ define
$$ \D x=\sum_{k\in\Z}kx_k.$$ Then $\D:X_\D\to X$ is a
is self-adjoint, regular operator on $X$.
\end{prop}

{\bf Remark} On the generators of the graph $C^*$-algebra we have the formula
 $$\D(S_\al S_\beta^*)=(|\al|-|\beta|)S_\al S_\beta^*.$$

The operator $\D$ gives us an unbounded Kasparov module, but we work with the
bounded Kasparov module we obtain from the phase of $\D$. We need a
preparatory
Lemma.

\begin{lemma}\label{finrank} Assume that the directed graph $E$ is
locally finite and has no sources. For all $a\in A$ and $k\in\Z$,
$a\Phi_k\in End^0_F(X)$, the
 compact right endomorphisms of $X$. If $a\in A_c$ then
$a\Phi_k$ is finite rank.
\end{lemma}

In fact we show that for $k\geq 0$, with $|v|_k=$the number of paths of length $k$ ending at $v\in E^0$,
$$\Phi_k=\sum_{|\mu|=k}\Theta_{S_\mu,S_\mu},\ \ \ \ \Phi_{-k}=\frac{1}{|v|_{k}}\sum_{|\mu|=k}\Theta_{S_\mu^*,S_\mu^*},$$
 where the sums converge strictly in $End_F(X)$ and $\Theta_{x,y}z=x(y|z)_R$.
In the particular case of the Cuntz algebra, which we examine later,
we have 
\be\Phi_{-k}=\frac{1}{n^k}\sum_{|\mu|=k}\Theta_{S_\mu^*,S_\mu^*}.\label{cuntzprojs}\ee

\begin{thm} Suppose that the graph $E$ is locally finite and has no sources,
and let $X$ be the right $F$ module of Definition \ref{mod}. Let
$V=\D(1+\D^2)^{-1/2}$. Then $(X,V)$ is an odd Kasparov module for
$A$-$F$ and so defines an element of $KK^1(A,F)$.
\end{thm}

The next theorem presents a general result about the Kasparov product
in the odd case.

\begin{thm}\label{themapH} Let $(Y,T)$ be an odd Kasparov module for the
$C^*$-algebras $A,B$. Then the Kasparov product of $K_1(A)$ with the
class of $(Y,T)$ is represented by
$$\la [u],[(Y,T)]\ra= [\ker PuP]-[{\rm coker} PuP]\in K_0(B),$$
 where $P$ is the non-negative spectral projection for $T$.
\end{thm}

Note that the above construction and Theorems 4.6 and 4.7 do not require
a trace.We also remark that
this pairing was studied in \cite{pr}, as well as the relation to the
semifinite index. To discuss all of this we must
first build a semifinite spectral triple: this is our next task.

\section{Semifinite spectral triples for graph algebras}

In order to obtain such a spectral triple, we require that there
exists on $A=C^*(E)$ a faithful, semifinite, norm
lower-semicontinuous, gauge invariant trace. In \cite{pr} we studied
the question of the existence of such traces, and showed that they were in
one-to-one correspondence with `graph traces' on $E$. These are
positive real-valued functions defined on the vertices of $E$ which
reflect the structure of the graph. The most notable of the various
necessary conditions for the existence of such a trace is that, in the
graph $E$, no loop should have an exit. See \cite{pr} for more
information.

We will begin with the right $F_c$ module $X_c$. In order to deal
with the spectral projections of $\D$ we will also assume
throughout this section that $E$ is locally finite and has no
sources. This ensures, by Lemma \ref{finrank}, that for all $a\in
A$ the endomorphisms $a\Phi_k$ of $X$ are compact endomorphisms. We
also assume that there exists a faithful, semifinite, norm
lower-semicontinuous, gauge invariant trace $\tau$ on $C^*(E)$
which by \cite{pr} is finite on $A_c$.

We
define a ${\C}$-valued inner product on $X_c$:
$$ \la x,y\ra:=\tau((x|y)_R)=\tau(\Phi(x^*y))=\tau(x^*y).$$
We define the
Hilbert space $\HH=L^2(X,\tau)$ to be the completion of $X_c$ for
$\la\cdot,\cdot\ra$. We need some additional information
about this situation in order to construct a 
spectral triple.  All of the following statements are proved in \cite{pr}.

First the $C^*$-algebra $A=C^*(E)$ acts on $\HH$
by an extension of left multiplication. This defines a faithful nondegenerate
$*$-representation of $A$. Moreover, any endomorphism of $X$ leaving $X_c$
invariant extends uniquely to a bounded linear operator on $\HH$.
Also the endomorphisms $\{\Phi_k\}_{k\in\Z}$ define mutually orthogonal projections on
$\HH$. For any $K\subset \Z$ the sum $\sum_{k\in K}\Phi_k$
converges strongly to a projection in $\B(\HH)$. In particular,
$\sum_{k\in\Z}\Phi_k=Id_{\HH}$, and so for all $x\in \HH$ the sum
$\sum_k\Phi_kx$ converges in norm to $x$.

Next the unbounded self adjoint operator needed for our spectral triple
is the 
operator $\D$ with domain $X_c$ extended to a
closed densely defined self-adjoint
operator on $\HH$.
Then the  technical details needed for a smooth algebra are handled
by the following results.

\begin{lemma}\label{deltacomms} Let $\HH,\D$ be as above and let
$\dd=\sqrt{\D^*\D}=\sqrt{\D^2}$ be the absolute value of $\D$.
Then for $S_\al S_\beta^*\in A_c$, the operator $[\dd,S_\al
S_\beta^*]$ is well-defined on $X_c$, and extends to a bounded
operator on $\HH$ with
$$\Vert[\dd,S_\al S_\beta^*]\Vert_{\infty}\leq
\Bigl||\al|-|\beta|\Bigr|.$$ Similarly, $\Vert[\D,S_\al
S_\beta^*]\Vert_\infty= \Bigl||\al|-|\beta|\Bigr|$.  \end{lemma}

\begin{cor}\label{smodense} The algebra $A_c$ is contained in the smooth domain
of the derivation $\delta$ where for $T\in\B(\HH)$, $\delta(T)=[\dd,T]$. That is $$ A_c\subseteq\bigcap_{n\geq 0}{\rm dom}\ \delta^n.$$
\end{cor}

\begin{defn} Define the $*$-algebra $\A\subset A$ to be the
completion of $A_c$ in the $\delta$-topology. By \cite[Lemma 16]{R1},
$\A$ is Fr\'{e}chet and stable under the holomorphic functional
calculus.
\end{defn}
The next Lemma addresses key facts necessary for dealing with
summability issues of spectral triples for nonunital algebras. We will
not discuss this here, but merely point out that the Lemma allows us
to apply the results of \cite{R2}. For more details on the application
in this setting see \cite{pr}.
\begin{lemma}\label{smoalg} If $a\in\A$ then $[\D,a]\in\A$ and the operators $\delta^k(a)$,
$\delta^k([\D,a])$ are bounded for all $k\geq 0$. If $\phi\in F\subset\A$ and
$a\in\A$
satisfy $\phi a=a=a\phi$, then $\phi[\D,a]=[\D,a]=[\D,a]\phi$.
The norm closed algebra generated by $\A$ and $[\D,\A]$ is $A$. In
particular, $\A$ is quasi-local.
\end{lemma}

We now come to the two key definitions of this Section. 
\begin{defn} Let $End^{00}_F(X_c)$ denote the algebra of finite rank
operators on $X_c$ acting on $\HH$. Define
$\cn=(End^{00}_F(X_c))''$, and let $\cn_+$ denote the positive
cone in $\cn$.
\end{defn}

We need an appropriate trace on $\mathcal N$. Our object is to define
one that is naturally related to $\tau$.

\begin{defn} Let  $T\in\cn$ and $\mu\in E^*$. Let $|v|_k=$ the
 number of paths of
 length $k$ with range $v$, and define for $|\mu|\neq 0$
$$\omega_\mu(T)=
\la S_\mu,TS_\mu\ra+\frac{1}{|r(\mu)|_{|\mu|}}\la
S_\mu^*,TS_\mu^*\ra.$$ For $|\mu|=0$, $S_\mu=p_v$, for some $v\in
E^0$, set $$\omega_\mu(T)=\la S_\mu,TS_\mu\ra=
\la p_v,Tp_v\ra =
\tau(p_vTp_v).$$
 Define
$$\tilde\tau:\cn_+\to[0,\infty],\ \ \mbox{by}\ \ \ \ \tilde\tau(T)=
\lim_{L\nearrow}\sum_{\mu\in L\subset E^*}\omega_\mu(T)$$
where $L$ is in the net of finite subsets of $E^*$.
\end{defn}

The immediate consequence of this definition
is that the functional $\tilde\tau$ satisfies: for $T,S\in\cn_+$ and 
$\lambda\in\R$, $\lambda\geq 0$ 
$$\tilde\tau(T+S)=\tilde\tau(T)+\tilde\tau(S)\ \ \
\mbox{and}\ \ \ \tilde\tau(\lambda T)=\lambda\tilde\tau(T)$$
where we adopt the convention $0\cdot\infty=0$.
Now we have the main property of $\tilde\tau$.
\begin{prop}\label{tildetau}
 The functional $\tilde\tau:\cn_+\to[0,\infty]$ defines a faithful normal
semifinite trace on $\cn$. Moreover,
$$End_F^{00}(X_c)\subset\cn_{\tilde\tau}:=
{\rm span}\{T\in\cn_+:\tilde\tau(T)<\infty\},$$
the domain of definition of $\tilde\tau$, and
$$\tilde\tau(\Theta^R_{x,y})=\la y,x\ra=\tau(y^*x),\ \ \ x,y\in X_c.$$
\end{prop}

We are now in a position to see that we have the semifinite
spectral triple that we promised.

\begin{thm}\label{mainthm} Let $E$ be a locally finite graph
with no sources, and let $\tau$ be a faithful, semifinite, gauge
invariant, lower semicontinuous trace on $C^*(E)$. Then
$(\A,\HH,\D)$ is a $QC^\infty$, $(1,\infty)$-summable, odd, local,
semifinite spectral triple (relative to $(\cn,\tilde\tau)$). For
all non-zero $a\in \A$, the operator $a(1+\D^2)^{-1/2}$ is not trace class.
If $v\in E^0$ has no sinks downstream
$$\tilde\tau_\omega(p_v(1+\D^2)^{-1/2})=2\tau(p_v),$$
where $\tilde\tau_\omega$ is any Dixmier trace associated to
$\tilde\tau$.
\end{thm}

Finally we can compute our numerical index pairing.

\begin{prop}\label{C*specflow} Let $E$ be a locally finite graph with
no sources and a faithful graph trace $g$, and $\tau_g$ the
corresponding faithful, semifinite, norm lower semicontinuous, gauge
invariant trace on 
 $A=C^*(E)$.  The pairing between the spectral triple
$(\A,\HH,\D)$ of Theorem \ref{mainthm} with $K_1(A)$ can be
computed as follows. Let $P$ be the positive spectral projection
for $\D$, and perform the $C^*$ index pairing of Theorem
\ref{themapH}:
$$K_1(A)\times KK^1(A,F)\to K_0(F),\ \ \ \
[u]\times[(X,P)]\to [\ker PuP]-[{\rm coker}PuP].$$ Then we have
$$sf(\D,u\D u^*)=\tilde\tau_g(\ker PuP)-\tilde\tau_g({\rm coker}PuP)=
\tilde\tau_{g}([\ker PuP]-[{\rm
 coker}PuP]).$$
\end{prop}

In fact this result is a special case of a much more general
relationship which we describe in the next section.

\section{The General Relationship Between Semifinite Index Theory and 
 $KK$-Theory}

A finitely summable spectral triple $(\A,\HH,\D)$ (in the usual
$\B(\HH)$ sense) represents a $K$-homology class $[(\A,\HH,\D)]\in
K^*(\overline{\A})$ for the $C^*$-algebra
$\overline{\A}$, \cite{C,CP1}. The local index theorem in this case
computes the pairing between this class and the $K$-theory 
$K_*(\overline{\A})$ in terms of the pairing of cyclic homology and cohomology
via Chern characters. 

For a semifinite spectral triple $(\A,\HH,\D)$ defined relative to
$(\cn,\tau)$,  
it turns out that there is always  a $KK$-class $[(\A,\HH,\D)]\in
KK^*(\overline{\A},B)$, where $B$ is constructed from the data of the
spectral triple. Thus there is a $K_*(B)$-valued index given by the
Kasparov product with the $K$-theory of $\overline{\A}$:
$$K_*(\overline{\A})\times KK^*(\overline{\A},B)\to K_*(B).$$
If the spectral triple is odd (even) and we pair with odd (even)
$K$-theory of $\overline{\A}$ we obtain an index in $K_0(B)$. Let us
just consider the odd-odd pairing for a moment. Let $P$ be the
non-negative spectral projection of $\D$, $u\in\A$ unitary, and denote by
$\mbox{Index}(PuP)$ the class in $K_0(B)$ given by Theorem
\ref{themapH}. Then it turns out that $PuP$ defines a Fredholm operator in the von
Neumann algebra $P\cn P$, and so we may take the trace of its kernel and
cokernel projections. Then the analytic pairing yields a number, the spectral
flow, and under some additional hypotheses, it is given by 
$$sf_\tau(\D,u\D
u^*)=\tau_*(\mbox{Index}(PuP))=\tau(\ker(PuP))-\tau(\mbox{coker}(PuP)).$$
This statement holds under the condition that certain residue traces
on $\A$ are faithful.
An analogous statement holds in the even case. These results are
proved in \cite{KNR}.

We may translate all this to cyclic theory via Chern characters:

$$\begin{array}{ccccc}
K_1(\overline{\A})&\times&KK^1(\overline{\A},B)&\to&K_0(B)\\
\downarrow\ Ch & &\downarrow\ Ch & &\downarrow\ Ch\\
H_1(\A)&\times & H^1(\A,\B)&\to&H_0(\B)\end{array}$$
Since the Chern characters respect the pairing, we may pair in
$K$-theory and apply the Chern character, or vice versa.
If, as is the case in the applications, we may regard the kernel and
cokernel projections of $PuP$ as living in an algebra $\K(B)$ Morita
equivalent to $B$, and with $\K(B)\subset\cn$, then we may regard the
trace $\tau$ as a cyclic zero cocycle on $\K(B)$, and compute the
pairing. Again, we compute the analytic spectral flow:
\bean\la Ch(u)\times Ch(\A,\HH,\D),\tau\ra&=&\la Ch(u\times
(\A,\HH,\D)),\tau\ra\nno
& =&\tau(\ker(PuP))-\tau(\mbox{coker}(PuP))\nno
&=&sf_\tau(\D,u\D
u^*).\eean
Thus the semifinite index theorem is computing a `zero-order' part of
a more general index pairing in bivariant cyclic theory, where $\tau$
may be replaced by a more general cocycle.


\section{Graph Algebras as Noncommutative Manifolds}

In \cite{PRen} we consider natural generalisations of Connes' axioms for
noncommutative manifolds, \cite{C1}, which simultaneously address semifinite
triples and nonunital algebras. The axiom set we use is modelled on
that of \cite{RV} where the spin manifold theorem is proven.

Most of the axioms admit straightforward generalisations, and we will
not dwell on them here, other than to say that most are automatically
satisfied by the gauge spectral triple of a locally finite graph $E$ with
no sources and such that $C^*(E)$ carries a faithful, semifinite, norm
lower semi-continuous, gauge invariant trace. In particular, the graph
$E$ cannot have a loop with an exit, \cite{pr}. The most interesting
axioms are those which impose
further constraints on the graph.

The first, and most restrictive, axiom we consider is
orientability. This axiom requires the existence of a Hochschild cycle $c$ of the same degree as the (integral) dimension, which in this
case is 1, such that $\pi_\D(c)=1$ (see below). We define
$$c=\sum_{e\in E^1}S_e^*\otimes S_e.$$
 Here the sum is over all edges in the graph. If $E$ is the graph with
a single vertex and single edge, then $C^*(E)\cong C(\T^1)$ and $c$ is
the usual volume form. In general, the convergence of the above sum
(in the tensor product of the multiplier algebra with itself) and the
condition $b(c)=0$ are both guaranteed by the following `single entry
hypothesis' 
$$\mbox{every vertex receives precisely one edge, and emits at least one.}$$
Coupled with the fact that no loop may have an exit, we see that the
graphs we may consider are (disjoint unions) of two types:
$$\mbox{infinite directed trees and a single loop comprised of}\ N\mbox{ vertices and edges}$$
 The latter type are isomorphic to $M_N(C(\T^1))$ while the former are
always nonunital and AF.

The Hochschild cycle $c$ may be interpreted in several ways, but the
most appealing to us is as a limit of compactly supported Hochschild
1-cycles, as in \cite{R1}. This is what one would expect of the (class
of) a volume form in de Rham theory.

The representation of $c$, $\pi_\D(c)$, is by definition
$$\pi_\D(c)=\sum_{e\in E^1}S_e^*[\D,S_e]=\sum_{e\in E^1}p_{r(e)}$$
where $S^*[\D,S_e]=p_{r(e)}$ is a straightforward calculation.
Since for each vertex $v$ there is a unique edge $e$ with range $v$, we
can show that the sum on the right converges strictly to the identity
in the multiplier algebra of $C^*(E)$, and strongly to the identity on
$\HH$. Thus the orientability condition is satisfied for directed
trees and finite loops.

The next axiom to consider is finiteness. Essentially this asks that the
smooth domain of $\D$, $\HH_\infty=\cap_{m\geq 1}\mbox{dom}\D^m$, be a
finite projective module over $\A$. Thus $\HH_\infty$ should be
thought of as a module of smooth sections of a `noncommutative vector
bundle'. 

Semifiniteness does not affect us
here, but nonunitality does; see \cite{R1}. In particular, there are
no issues to consider here for the $N$-point loop examples, only the
directed trees. We take the approach that
we should be able to recover the `continuous sections vanishing at
infinity' from these `smooth sections'.

We summarise what our axiom of finiteness says:

1) $\HH_\infty$ should be a continuous $\A$-module, \vspace{-6pt}

2) $\HH_\infty$ embeds continuously as a dense subspace in the
$C^*$-module $X$,\vspace{-6pt}

3) $X$ is the completion of $pA^N$ for some $N$
($A^N=A\oplus\ldots \oplus A$, $N$ copies) and some
projection $p$ in the $N\times N$ matrices over a unitization of $A=C^\ast(E)$, 
\vspace{-6pt}

4) the Hermitian product $\HH_\infty\ni x,y\to x^*y$ should have
range in $\A$ (acting on the right).

 This gives a relationship
between `decay at infinity' in the $C^*$-module sense and the Hilbert
space sense. In particular, the trace must be bounded below, which in
turn means that the graph can have only finitely many branchings. This
is then intimately related to the $K$-theory of the graph algebra.

\begin{prop}\label{finprojcore}  Suppose that the locally
finite directed graph $E$ has no sinks, no loops and satisfies the
single entry condition.
The $\A$-module $\HH_\infty$ satisfies 1) and 3). The module
$\HH_\infty$ satisfies  2)  if and only if the
$K$-theory of $A=C^\ast(E)$ is finitely generated. In this case the Hilbert
space $\HH$ also satisfies the conditions 2). If the $K_*(A)$ 
is finitely generated then condition 4) holds.
\end{prop}

Thus when we consider the directed tree examples, we must restrict to
those with finitely many branchings. It is interesting to observe that
the name `finiteness' for this axiom really refers to the finite
projective module. It is thus surprising to see this axiom influence
also the finiteness of the rank of the $K$-theory. 

The final axiom we mention is Poincar\'{e} Duality. While it was shown
in \cite{RV} that this axiom could be replaced by `closedness' and a
`spin$^c$' condition (versions of both are satisfied by the gauge spectral triples of graph algebras), Poincar\'{e} Duality is still an
important structural feature to examine. 

It turns out that the algebra $A=C^\ast(E)$ 
does satisfy Poincar\'{e} Duality in
 $K$-theory, though one needs a suitable nonunital formulation,
\cite{R1}. However, more is true. The fixed point subalgebra $F$ for
the gauge action also satisfies Poincar\'{e} Duality. This is true
both for the $N$-point loops and the directed trees (with finitely
many branchings). In essence, Poincar\'{e} Duality for {\em both}
these algebras arises as a consequence of the other axioms.

In \cite{PRen} we also consider the examples of higher rank graphs, or
$k$-graphs, studied from this point of view in \cite{prs}, where
$(k,\infty)$-summable semifinite spectral triples were obtained under
similar hypotheses on the $k$-graph as in the graph case. Again, most
of the axioms are satisfied immediately, with orientation and
finiteness placing similar constraints on the algebra (or
 $k$-graph). We have not considered Poincar\'{e} Duality for these
higher rank graphs, but expect phenomena similar to that for the
1-graph case.

\section{The Cuntz Algebra}\label{cuntzsection}

\subsection{The KMS state for the gauge action}

The Cuntz algebra does not possess a faithful gauge invariant
trace. It does however have a unique KMS state 
relative to the gauge action which is given by
the faithful gauge invariant state $\tau\circ\Phi:O_n\to\C$
where $\Phi:O_n\to F$ is the expectation and $\tau:F\to\C$ the unique
faithful normalised trace on the gauge invariant subalgebra.

Since the Cuntz algebra is the graph algebra of a locally finite
graph with no sources, the generator of the gauge action $\D$
acting on the right $C^*$-$F$-module $X$ gives us a Kasparov module
$(X,\D)$. As with tracial graph algebras, we take
this class as our starting point.

The first difficulty we encounter is that there are no unitaries to
pair with, since $K_1(O_n)=0$.
We require a new approach.

Let $\HH=L^2(O_n)$ where the inner product is defined by
$$\la a,b\ra=(\tau\circ\Phi)(a^*b).$$
Then $\D$ extends to a self-adjoint unbounded operator on $\HH$, \cite{PRen},
and we denote this closure by $\D$ from now on. The representation
of $O_n$ on $\HH$ (by left multiplication) is bounded and
nondegenerate, and the dense subalgebra $\mbox{span}\{S_\mu
S_\nu^*\}$ is in the smooth domain of the derivation $\delta$.
Thus we
see that the central algebraic structures of the gauge spectral
triple on a tracial graph algebra are mirrored in this
construction.

What differs significantly from the tracial situation is the
analytic information. We begin by obtaining some information about
the trace on $F$, the corresponding state on $O_n$ and the
associated modular theory.

First of all the  trace $\tau:F\to\C$ satisfies 
\ben \tau(S_\mu S^*_\nu)=\delta_{\mu\nu}\frac{1}{n^{|\mu|}}.\een
In terms of the associated state 
on $ O_{nc}$ we have some additional structure:

\begin{lemma} The  algebra $ O_{nc}={\rm span}\{S_\mu S_\nu^*\}$ with
the inner product arising from the state $\tau\circ\Phi$ is a modular
Hilbert algebra. \end{lemma}

Let $S$ first denote the operator
$a \mapsto a^*$ defined on $O_{nc}$. 
The adjoint $F=S^*$ of $S$ can be explicitly calculated on $O_{nc}$
and satisfies:
 $$F(S_{\mu}S_{\nu}^*)=n^{(|\mu|-|\nu|)}S_{\nu}S_{\mu}^*.$$
In particular, $F$ is densely defined so that $S$ is closable and we use
$S$ to denote the closure of $S$ restricted to $O_{nc},$
and also $F$ will denote the closure of $F$ restricted to $O_{nc}$.
 Then $S$  has a polar decomposition as 
\ben S=J\Delta^{1/2},\ \  \Delta=S^*S,\een 
where $J$
is an antilinear map, $J^2=1$. The Tomita-Takesaki modular theory, \cite{KR},
shows that \ben \Delta^{-it}O_n''\Delta^{it}=O_n'',\ \ \
JO_n''J^*=(O_n'')',\een where $O_n''$ is the weak closure of
$O_n$ acting on $L^2(O_n,\tau\circ\Phi)$.

\begin{lemma} Writing $S:a\to a^*$ as $S=J\Delta^{1/2}$, we have
\ben \Delta S_\mu S^*_\nu=n^{|\nu|-|\mu|}S_\mu S^*_\nu,\ \ \ \
\forall\mu,\nu\in{\N}^n.\een Thus the group of modular
automorphisms is given by \ben \s_t(S_\mu
S^*_\nu)=\Delta^{it}S_\mu
S^*_\nu\Delta^{-it}=n^{it(|\nu|-|\mu|)}S_\mu S^*_\nu.\een
\end{lemma}

\begin{cor} If $\D$ is the generator of the gauge action of $\T^1$ on
$O_n$, 
we have
\ben \Delta=e^{-\D\log n}\  {\rm  or  }\
e^{it\D}=e^{-\frac{it}{\log n}\log\Delta}.\een \end{cor}

Hence the flows generated by $\log\Delta$ and $\D$ are the same, up
to a constant rescaling of $-\log n$.  More precisely
\ben \Delta^{it}=e^{it\log\Delta}=e^{-(it\log n)\D}.\een
To continue, we recall the underlying right $C^*$-$F$-module, $X$,
which is the completion of $O_n$ for the norm $\Vert
x\Vert_X^2=\Vert\Phi(x^*x)\Vert_F$. Endomorphisms of the module $X$
which preserve $X_c$ (the copy of $O_{nc}$ inside $X$) extend uniquely
to bounded operators on the Hilbert space $\HH$, \cite{pr}. 

\subsection{An attempted semifinite spectral triple}

We follow our previous strategy (Section 5) to uncover the difficulty.

\begin{prop}\label{tildetrace} Let $\cn$ be the von Neumann algebra
$$  \cn=(End^{00}_F(X_c))'',$$
where we take the commutant inside $\B(\HH)$. Then $\cn$ is 
semifinite and there exists a faithful,
semifinite, normal trace $\tilde\tau:\cn\to\C$ such that for all rank
one endomorphisms $\Theta_{x,y}$ of $X_c$, 
$$\tilde\tau(\Theta_{x,y})=(\tau\circ\Phi)(y^*x),\ \ \ x,y\in X_c.$$
 In addition, $\D$ is affiliated to $\cn$.
\end{prop}

Rather than explain the proof here, we simply observe for the reader's
benefit that to check the trace
property (on endomorphisms) only requires that $\tau$ is a trace on 
$F$, not all of $O_n$. Here is the formal
calculation for rank one operators.
\bean \tilde\tau(\Theta_{w,z}\Theta_{x,y})&=&\tilde\tau(\Theta_{w(z|x),y})
=\tau((y|w(z|x)))\nno
&=&\tau((y|w)(z|x))=\tau((z|x)(y|w))\nno
&=&\tilde\tau(\Theta_{x(y|w),z})=\tilde\tau(\Theta_{x,y}\Theta_{w,z}).\eean

However this trace is not what we
need for defining summability.
To see this  we do some calculations. For $k\geq 0$ 
$$
\tilde\tau(\Phi_k)=\tilde\tau(\sum_{|\rho|=k}\Theta_{S_\rho,S_\rho})
=\tau(\sum_{|\rho|=k}(S_\rho|S_\rho))=\tau(\sum_{|\rho|=k}S_\rho^*S_\rho)=\sum_{|\rho|=k}1
=n^k.$$ Similarly, for $k<0$
$$\tilde\tau(\Phi_k)=n^k.$$
With respect to this trace it is not hard to see 
that we cannot expect $\mathcal D$ to satisfy
any summability criterion. To get summability we
are forced to define a new functional on $\cn$ by
$$\tau_\Delta(T):=\tilde\tau(\Delta T).$$
Since $\tilde\tau$ is a faithful semifinite normal trace on $\cn$, and
$\Delta$ is a positive invertible operator affiliated to $\cn$
we may show that
$\tau_\Delta$ is a faithful semifinite normal weight on $\mathcal N$.

To understand what we have done here we do
some further computations: for $k\geq 0$
$$
\tau_\Delta(\Phi_k)=\tau(\sum_{|\rho|=k}(S_\rho|\Delta S_\rho))
=\sum_{|\rho|=k}n^{-k}\tau(S^*_\rho S_\rho)=1\ \ \mbox{and}\ \ 
\tau_\Delta(\Phi_{-k})=1.$$
However, $\tau_\Delta$ is not a trace on finite rank endomorphisms. \bean
\tau_\Delta(\Theta_{x,y}\Theta_{w,z})&=&\tau_\Delta(\Theta_{x(y|w),z})
=\tilde\tau(\Delta\Theta_{x(y|w),z})=\tau((z|\Delta x(y|w)))\nno
&=&\tau((z|\Delta x)(y|w))\ \qquad\qquad\mbox{since}\ \Delta\
\mbox{is linear over}\ F\nno
&=&\tau((y|w)(z|\Delta x))\ \qquad\qquad\mbox{since}\ \tau\ \mbox{is a
trace on}\ F\nno
&=&\tau((y|w)(\Delta z|x))\ \qquad\qquad\mbox{since}\ \Delta\
\mbox{is self-adjoint on }X\nno
&=&\tilde\tau(\Theta_{w,\Delta z}\Theta_{x,y})\nno
&=&\tilde\tau(\Delta\Delta^{-1}\Theta_{w,z}\Delta\Theta_{x,y})
=\tau_\Delta(\Delta^{-1}\Theta_{w,z}\Delta\Theta_{x,y}).\eean

\begin{lemma} The modular automorphism group $\s_t^{\tau_\Delta}$ of
$\tau_\Delta$ is inner and given by
$\s_t^{\tau_\Delta}(T)=\Delta^{it}T\Delta^{-it}$. The weight
$\tau_\Delta$ is a KMS weight for the group $\s_t^{\tau_\Delta}$, and 
$$\s_t^{\tau_\Delta}|_{O_n}=\s_t^{\tau\circ\Phi}.$$
\end{lemma}


We cannot construct a semifinite spectral triple for 
$\mathcal N$. The best we can do is work with
 $\cM=\cn^{\s}$ the fixed point algebra for the modular automorphism group of $\tau_\Delta$.
Then any  $m\in\cM$ commutes with the spectral projections of $\Delta$ or 
\ben\label{commutes}  [\Delta,m]=0.
\een
Also observe that since $F$ and the projections $\Phi_k$ are invariant,
 they belong to $\cM$.
Finally the weight $\tau_\Delta$ restricts to a normal,
faithful, semifinite trace on $\cM$.

The reward for having sacrificed a trace on $\cn$ for a trace on $\cM$
is the following.

\begin{prop}\label{dixycomp} We have 
$(1+\D^2)^{-1/2}\in\LL^{(1,\infty)}(\cM,\tau_\Delta)$, and for
all $f\in F$ and all Dixmier traces $\tau_{\Delta,\omega}$
$$ \tau_{\Delta,\omega}(f(1+\D^2)^{-1/2})=2\tau(f).$$
\end{prop}

\begin{cor}
For all $a,b\in O_{nc}$,
$$lim_{s\to 1}(s-1)\tau_\Delta(ab(1+\D^2)^{-s/2})=2\tau(ab)=
2\tau(\s_i(b)a)=
\tau_{\Delta,\omega}(\s_i(b)a(1+\D^2)^{-1/2}).$$
Moreover, the functional $a\to \tau_{\Delta,\omega}(a(1+\D^2)^{-1/2})$
is a KMS$_1$  state on $O_{nc}$.
\end{cor}

Before continuing we present a set of data which appear to be essential
 if we are to describe in general the
situation presented here for the Cuntz algebras; compare \cite[Definition 2.1]{Gos}.

A modular spectral triple $(\A,\HH,\D)$, relative to a semifinite von
Neumann algebra $\cn$ and faithful KMS state 
$\tau$ on $\A$ with respect to an action $\s$ of the circle, 
would appear to need the following data.

The  algebra $\cn$ acts on the GNS Hilbert space
$\HH$ associated to $\tau$ and

0) the $*$-algebra $\A$ is faithfully represented in $\cn$,

1) there is a faithful normal semifinite weight $\phi$
on $\cn$ such that the modular automorphism group of $\phi$ is an
inner automorphism group $\tilde\sigma$ of $\cn$ with
$\tilde\sigma|_\A=\s$, 

2) $\phi$ restricts to a faithful  semifinite trace on  $\cM=\cn^\s$,

3) $[\D,a]$ extends to a bounded operator (in $\cn$) for all $a\in\A$
and for $\lambda$ in the resolvent set of $\D$,
$f(\lambda-\D)^{-1} \in {\mathcal K}(\cM,\phi)$, where  $f\in\A^\s$, and
${\mathcal K}(\cM,\phi)$ is the ideal of compact operators in $\cM$ relative to
$\phi$. In particular, $\D$ is affiliated to $\cM$.

The triple is even if there exists $\gamma=\gamma^*$, $\gamma^2=1$
such that $\gamma\D+\D\gamma=0$ and $\gamma a=a\gamma$ for all
$a\in\A$. Otherwise it is odd.

Observe that this data implies that $\D$ commutes with 
the modular automorphism group in the strong sense that its
spectral projections are invariant.

We will refer to any triple  $(\A,\HH,\D)$ satisfying 
this data as a modular spectral triple
recognising that as more examples are studied modifications 
or extensions may be required. In particular, this definition does not
address actions of the reals that do not factor through the
circle. 

Of course the triple $(O_{nc},\HH,\D)$ along with $\cn,\tau_\Delta$ 
constructed in this Section, is a modular spectral triple.

%




\subsection{Modular unitaries}

As the Cuntz algebra is not contained in $\mathcal M$ we cannot
use the semifinite spectral triple immediately. We need some
additional structure.

\begin{defn} Let $\A$ be a $*$-algebra and $\s:\A\to \A$ an algebra
automorphism such that 
$$ \s(a)^*=\s^{-1}(a^*).$$
then we say that $\s$ is a regular automorphism, \cite{KMT}.
\end{defn}

Clearly for a modular
spectral triple the  automorphism $\s(a)=\Delta^{-1}a\Delta$ is regular.

\begin{defn} Let $u$ be a unitary in a matrix algebra over $\A$, and $\s:\A\to \A$ a regular
automorphism with fixed point algebra $F=\A^\s$. We say that $u$
satisfies the {\bf modular condition} with respect to $\s$ 
if both the operators
$$ u\s(u^*),\ \ \ u^*\s(u)$$
are in (a matrix algebra over) the algebra $F$. We denote by
$U_\s$ the set of modular unitaries obtained by taking the union over all matrix
algebras over $\A$..
\end{defn}

We are of course thinking of the case $\s(a)=\Delta^{-1}
a\Delta$, where $\Delta$ is the modular operator for some weight on
$A$. Hence the use of the terminology `modular unitaries'.
For unitaries in matrix algebras over $\A$ we use the
regular automorphism 
$\s\otimes Id_n$ to state the modular condition, where $Id_n$
is the identity of $M_n(\C)$.

\begin{defn} Let $u_t$ be a continuous path of modular unitaries such
that $u_t\s(u_t^*)$ and $u^*_t\s(u_t)$ are also continuous paths in
$F$. Then we say that $u_t$ is a modular homotopy, and say that $u_0$ and
$u_1$ are modular homotopic.
\end{defn}

\begin{lemma} The binary operation on modular homotopy classes in
$U_\s$
$$[u]+[v]:=[u\oplus v]$$
is abelian.
\end{lemma}

We can now also see why the usual proof that the inverse of $u$ is
$u^*$ in $K_1(A)$ is not available to us. This usual proof is as
follows. Observe that $u\oplus v=(u\oplus 1)(1\oplus v)\sim
(1\oplus u)(1\oplus v)=(1\oplus uv)$. Then we see that addition in
$K_1$ arises from multiplication of unitaries, and so
$[u]+[u^*]=[uu^*]=[1]=0$. However, while the homotopy from
$u\oplus 1$ to $1\oplus u$ is a modular homotopy 
in $U_\s$ by the last Lemma, the homotopy from $(u\oplus
1)(1\oplus v)$ to $(1\oplus u)(1\oplus v)$ is not in general. The
multiplication on the right by $(1\oplus v)$ breaks the
modular condition. In particular, the product of two modular unitaries
need not be a modular unitary.  What we do have in this situation is
stated in the following result.

\begin{lemma} If $u\in M_l(F)$ is unitary then $u\oplus u^*\sim
1$. If $u\in M_l(F)$ then $-[u]=[u^*]$ in $K_1(A,\s)$, 
with $K_1(A,\sigma)$ defined below.
\end{lemma}

We now formalise the above discussion. Compare the following
 with \cite[Definition 4.8.1]{HR}

\begin{defn} Let $K_1(A,\s)$ be the abelian semigroup with one
generator $[u]$ for each unitary $u\in M_l(A)$ satisfying
the modular condition
 and with the following relations:
\bean 1)&& [1]=0,\nno 2)&&
[u]+[v]=[u\oplus v],\nno 3)&& \mbox{If }u_t,\ t\in[0,1]\ \mbox{is a
continuous paths of unitaries in }M_l(A)\nno && \mbox{satisfying
the modular condition then}\ [u_0]=[u_1].\eean
\end{defn}

Of course we can make this into a group by the Grothendieck construction,
but as yet see no compelling reason to do so.

{\bf Example} For $S_\mu\in O_{nc}$ we write $P_\mu=S_\mu
S_\mu^*$. Then 
for each
$\mu,\nu$ we have a unitary

$$ u_{\mu,\nu}=\bma 1-P_\mu & S_\mu S_\nu^*\\ S_\nu S_\mu^* &
1-P_\nu\ema.$$ It is simple to check that this a self-adjoint
unitary satisfying the modular condition.


\begin{lemma}\label{noinv} For all $\mu,\nu$ there is a modular homotopy
$$ u_{\mu,\nu}\sim u_{\nu,\mu}.$$
\end{lemma}

{\bf Example} More generally, if $\s$ is a regular automorphism of an
algebra $A$ with fixed point algebra $F$, $v\in A$ is a partial
isometry with range and source projections in $F$, and furthermore
$v\s(v^*),\ v^*\s(v)\in F$, then 
$$u_v=\bma 1-v^*v & v^*\\ v & 1-vv^*\ema$$
is a modular unitary, as the reader may check. The proof
of Lemma \ref{noinv} applies to these unitaries to show that $u_v\sim
u_{v^*}$. 

There are a number of structural results we can prove about this situation.
For example it is not hard to see that the centre of $\cM$, denoted
${\mathcal Z}(\cM)$, contains the von Neumann algebra generated by the 
spectral projections of $\Delta$. A more crucial matter is worth
emphasising in the next result.



\begin{lemma} Let $(\A,\HH,\D)$ be a modular spectral triple such that
$\D$ commutes with $F=\A^\s$. Let $u\in\A$ be a
 unitary. Then $uQu^*\in\cM$ for all spectral projection $Q$ of $\D$,
 if and only if $u$ is modular. 
\end{lemma}

\begin{proof} First, $uQu^*$ is a projection in $\cn$. For one
direction  we
have
\bean \s(uQu^*)&=&\s(u)Q\s(u^*),\ \ \ \ \ \ \ Q\in\cM\nno
&=& uu^*\s(u)Q\s(u^*)\nno
&=&uQu^*\s(u)\s(u^*),\ \ \ \ \ \ \ u^*\s(u)\in F\nno
&=&uQu^*.\eean
Hence $uQu^*$ is invariant, and so in $\cM$. On the other hand if 
$$ uQu^*=\s(uQu^*)=\s(u)Q\s(u^*)$$
then we have 
$$Q=u^*\s(u)Q\s(u^*)u=Q+[u^*\s(u),Q]\s(u^*)u.$$
As $\s(u^*)u$ is invertible, we see that
$[u^*\s(u),Q]=0$. Since $u^*\s(u)\in\A$, and commutes with all $Q$, it
lies in $F=\cM\cap\A$.
\end{proof}

The truly important aspect of this lemma is that modular
unitaries conjugate $\Delta$ to an element of $\cM$, and so $u\Delta
u^*$ commutes with $\Delta$. 




\section{The Local Index formula for the
Cuntz Algebras}
To obtain a pairing  between modular
$K_1$ and modular spectral triples, we are going to use
the spectral flow formula of \cite{CP2}. 
Before we can do this
effectively, we need to address an important detail. 
The modular unitaries do not lie in the algebra
$\mathcal M$. Consequently 
the formula in \cite{CP2} for spectral flow
actually contains two additional terms
which measure the spectral asymmetry of the end-points.
These eta-type correction terms  cancel when the end-points
are unitarily equivalent, however the argument demonstrating this
fails in the present situation because modular unitaries lie
outside $\mathcal M$.
It eventuates, by a different argument, that
for modular unitaries of the form $u=u_v$ discussed
above, and under the hypotheses of the theorem below,
 these correction terms do cancel out. 
For  $r>0$ this gives us
$$ sf(\D,u\D u^*)=
\frac{1}{C_{1/2+r}}\int_0^1\phi(u[\D,u^*](1+(\D+tu[\D,u^*])^2)^{-1/2-r})dt.$$
We are now in a position to apply 
the reasoning in 
\cite{CPRS2}, to obtain a residue 
formula to compute the index pairing between
$( O_{nc},\HH,\D)$ and $K_1(O_n,\s)$ for these particular
unitaries. Let $\D_k:=\D\otimes Id_k$ and $\phi_k=\phi\otimes \mbox{Tr}_k$.

\begin{thm}\label{first} Let
$( O_{nc},\HH,\D)$ be the $QC^\infty$, $(1,\infty)$-summable, 
modular spectral triple relative to $(\cn,\tau_\Delta)$, constructed previously.
 Then for any modular
unitary of the form $u_v$ with $u_v\in M_k(O_{nc})$, with $v$ a partial isometry with range and
source projections in $M_{k/2}(F)$ and $v\s(v^*),\ v^*\s(v)\in M_{k/2}(F)$, and any
Dixmier trace $\phi_{k,\omega}$ we have
$$sf_{\phi_k}(\D_k,u_v\D_k u_v)=lim_{r\to 0}
r\phi_k(u_v[\D_k,u_v](1+\D_k^2)^{-1/2-r})=\frac{1}{2}\phi_{k,\omega}(u_v[\D_k,u_v](1+\D_k^2)^{-1/2}).$$
The functional
$$ M_k(O_{nc})\otimes  M_k(O_{nc})\ni a_0\otimes a_1\to
lim_{r\to 0}r\phi_k(a_0[\D_k,a_1](1+\D_k^2)^{-1/2-r})$$
is a twisted $b,B$-cocycle. Moreover, the spectral flow 
depends only on the modular homotopy class
of $u_v$.
\end{thm}

Following the proof of the local index theorem in
\cite{CPRS2}, one can also obtain a twisted resolvent cocycle
although we will not discuss this here.

\begin{thm} Let $(O_{nc},\HH,\D)$ be the modular spectral triple of the
Cuntz algebra, and $u$ a modular unitary of the form
$u_{\mu,\nu}$. Then 
\bean sf_{\phi_2}(\D_2, u\D_2 u^*)&=&
(|\mu|-|\nu|)\left(\frac{1}{n^{|\nu|}}-\frac{1}{n^{|\mu|}}\right)
\in (n-1)\Z[1/n]\nno
&\geq&0.\eean
\end{thm}
Thus we see a wholly new kind of index pairing for the Cuntz
algebra. The index pairings above would be zero if we were employing a
trace, as all the unitaries we consider are self-adjoint, and so
represent zero in ordinary $K_1$. 

Modular $K_1$, and its pairing with twisted cyclic, represents a
new way to obtain data about algebras without (faithful) traces, by
using KMS states instead. The procedure we have outlined for the Cuntz
algebra is in fact quite general, and we have recently applied it to
the Haar state of $SU_q(2)$ with similar success, \cite{CRT}. 

We expect to find many applications of these new invariants, both in
mathematics and physics.



\end{document}